\documentclass[11pt]{article}
\usepackage[T1]{fontenc}
\usepackage[latin9]{inputenc}
\usepackage[a4paper]{geometry}
\usepackage{color}
\usepackage{amsmath}
\usepackage{setspace}
\usepackage{amssymb}

\usepackage{amscd} 
\usepackage{tikz-cd}

\usepackage[active]{srcltx}

\newcommand{\dparcial}[2]{ \frac{\partial #1}{\partial #2} }

\newcommand{\modif}[1]{#1}

\newcommand{\modifr}[1]{\textcolor{red}{}}

\def\nn{\nonumber}

\newtheorem{thm}{Theorem}[section]

\newtheorem{cor}{Corollary}[section]

\begin{document}


\title{\textbf{\Large Wei-Norman Factorization of Frenet-Serret equations and the Cardan Lock}}

\author{Julio Guerrero$^{a,b}$ and Rafael Florencio$^c$}


\maketitle
\begin{singlespace}

{\noindent \begin{center}
\emph{$^{a}$ 
Department of Mathematics, Faculty of  Experimental Sciences, University of Ja\'en, Campus Las Lagunillas s/n, 23071 Ja\'en, Spain}
\end{center}}

{\noindent \begin{center}
\emph{$^{b}$ 
Institute Carlos I of Theoretical and Computational Physics, University of  Granada, Fuentenueva s/n, 18071 Granada, Spain}
\end{center}}

{\noindent \begin{center}
\emph{$^{c}$ Department of Mathematics, Higher Politechnic School of Linares,  University of Jaen, Campus Cient\'{\i}fico Tecnol\'ogico de Linares,  Avda. de la Universidad, Ronda Sur, s/n, 23700 Linares, Ja\'en}
\end{center}}

\end{singlespace}

\begin{abstract}
Frenet-Serret equations for particles moving along differentiable curves, in the case of non-constant curvature and torsion, are considered in the framework of the Wei-Norman factorization of the propagator using the Lie algebra of the rotation group. It leads to a system of coupled linear equations that are highly coupled if Euler or Cardan angles are used, but that can be uncoupled if the Cartan basis in the Lie algebra is used, resulting essentially in a Ricatti equation.  However, since the rotations algebra is not a solvable algebra, the global Wei-Norman Theorem does not apply and singularities in the solutions of the system of equations appear at certain points, related to the \text{gimbal o Cardan lock}.  Notwidthstanding, a detailed study of the limit at these points shows that in most of the cases the propagator is well defined and gimbal lock is avoided. 
\end{abstract}
Keywords: Wei Norman factorization, Frenet-Serret equations.


\section{Introduction}
\label{sec:intro}

The Frenet--Serret equations \cite{Frenet,Serret} are among the
oldest and most widely used tools of differential geometry. They describe
how the orthonormal moving frame 
attached
to a regular, non-degenerate curve in $\mathbb{R}^3$ rotates as the curve is
traversed, the rate of rotation being governed by only two functions of the
arclength: the curvature $\kappa$ and the torsion $\tau$. By the
fundamental theorem of curves, $\kappa$ and $\tau$ determine the curve
up to a rigid motion, so integrating the Frenet--Serret equations is
equivalent to reconstructing a curve from its intrinsic invariants. Beyond
geometry, the same equations appear in the kinematics and control of
vehicles and articulated robots \cite{Seligbook,SeligNeedle,SeligRational}, where
the rotation group $SO(3)$ is enlarged to the Euclidean group $SE(3)$, and
they admit natural generalizations to other spaces.

When $\kappa$ and $\tau$ are constant the Frenet--Serret equations have
constant coefficients and are solved by a single matrix exponential, which
produces circular helices. As soon as $\kappa$ and $\tau$ are allowed
to vary, closed-form solutions are known only
in special situations, such as plane curves ($\tau=0$) or generalized
helices of Lancret type ($\tau/\kappa$ constant) \cite{libroAngel}. The
classical approach, due to Lie and Darboux, reduces the problem through a
complex change of variables to a Riccati equation followed by two
quadratures \cite{Struik}. Although the Riccati equation cannot be
integrated in general, it is one of the best understood nonlinear ordinary
differential equations, and this reduction remains the natural starting
point for any analytic treatment.

In this paper we revisit the Frenet--Serret equations from the point of view
of the Wei--Norman method \cite{WeiNorman1,WeiNorman2}. Writing the
equations as a linear initial value problem $Q'(s)=A(s)Q(s)$ on $SO(3)$, with
$A(s)\in\mathfrak{so}(3)$, the propagator $U(s,s_0)$ is sought as an ordered
product of exponentials of the elements of a basis of the Lie algebra. The
Wei--Norman theorem guarantees that this factorization exists in a
neighbourhood of the initial point and that the exponents satisfy a system
of nonlinear first-order equations determined solely by the coefficients of
$A(s)$ and the structure constants of the algebra. This reduces the original
problem, with nine unknowns, to a problem with three unknowns. The method has
a long history in the theory of Lie systems and the nonlinear superposition
principle \cite{Marmo,LieSystem}, and it has been
used recently to give a semiclassical interpretation of the singularities
that appear for the group $SU(1,1)$ \cite{WNGuerreroBerrondo}.

A central point of the paper is that, since $\mathfrak{so}(3)$ is not
solvable, the Wei--Norman factorization is only local. Each choice of
ordered basis defines a local chart of the second kind on the group
\cite{Varadarajan}, and the Wei--Norman equations become singular where
the Jacobian  of that chart vanishes. For $SO(3)$ parameterized
by Euler or Cardan angles \cite{Shuster,GeneralizedEuler}
these singularities are precisely the well-known \emph{gimbal lock}, while
for $SU(1,1)$ they correspond to caustics and are related to the Maslov
correction \cite{WNGuerreroBerrondo}. The singularities are a property of
the chart and not of the problem: the propagator of the Frenet--Serret
equations exists and is smooth for every $s$ whenever $\kappa$ and $\tau$
are continuous, and no single three-parameter chart can cover the whole of
$SO(3)$ \cite{Stuelpnagel}.

The paper is organized as follows. In Section~2 we review the Wei--Norman
factorization method, sketch the proof of the local theorem, and interpret
the singular points in terms of local charts on the Lie group. We also
introduce the notion of wrenches, borrowed from robotics, and show that
vanishing wrenches lead to structural constraints, independent of the
coefficients, which can be used for motion planning. In Section~3 we recall
the Frenet--Serret equations and their formulation as an initial value
problem on $SO(3)$, emphasizing that the Frenet--Serret frame has one null
wrench. Section~4 contains the main results. Using the real basis
$\{L_x,L_y,L_z\}$, we obtain the Wei--Norman equations for the
Cardan  and Euler  factorizations, characterize
their singular points, and compute the limit of the propagator at gimbal
lock, showing explicitly how one degree of freedom is lost. Using instead the
complex basis $\{L_+,L_-,L_x\}$, we recover the Lie--Darboux reduction: the
equations uncouple into a Riccati equation and two quadratures. We give the
reality constraints that the complex exponents must satisfy, the
linearization of the Riccati equation into a complex harmonic oscillator
with $s$-dependent frequency, and a real parameterization in which $g_\pm(s)$
are stereographic coordinates of points on the Riemann sphere. In this
picture the singularity of the factorization occurs when the unit tangent
satisfies $\vec T(s)=-\vec T(s_0)$.


\section{The Wei-Norman factorization method}
\label{WN}

Consider the first-order linear initial value problem (IVP), where the prime indicates derivative with respect to $t$:
\begin{equation}
 y'(t)=A(t)y(t)\,,\qquad y(t_0)=y_0 \label{PVIvec}
\end{equation}
with $y(t),y_0\in V$ and $y_0$ the initial condition at $t_0$. Here $V$ is a vector space (which can be infinite-dimensional) and $A(t)$ is a family of endomorphism of $V$ (i.e. a family of matrices or  linear operators).


By Picard's existence and uniqueness theorem \cite{Picard} (particularized to systems of homogeneous linear equations), if $A(t)$ is continuous and $||A(t)||$ is bounded on some interval ${\cal I}\subset\mathbb{R}$ containing $t_0$, then the IVP \eqref{PVIvec} has a unique solution on ${\cal I}$.


It is convenient to write the solution $y(t)$   as:
\begin{equation}
 y(t)=U(t,t_0)y(t_0)\,,\qquad\forall t\in {\cal I}
\end{equation}
where $U(t,t_0)$ is known as the \textit{matrizant} \cite{Baker}, evolution operator, propagator, transfer matrix, etc., depending on the context.

Using the evolution operator, Eq.  \eqref{PVIvec} is transformed into
 \begin{equation}
 \dparcial{U(t,t_0)}{t}={A}(t){U}(t,t_0)\,,\qquad {U}(t_0,t_0)={I}_V \label{PVIop}
 \end{equation}
with ${I}_V$ the identity automorphism of $V$.  For $V$ finite-dimensional, it is easily checked that $U(t,t_0)$ is an automorphism of $V$ for all $t,t_0\in{\cal I}$, since we have that (see, for instance, \cite{MagnusReport})
\begin{equation}
 \det U(t,t_0)=\det U(t_0,t_0)\exp\left(\int_{t_0}^t {\rm tr} (A(s))ds\right) =\exp\left(\int_{t_0}^t {\rm tr} (A(s))ds\right)
 \label{detU}
\end{equation}
since $\det U(t_0,t_0)=1$ and $A(t)$ is bounded on ${\cal I}$. The infinite dimensional case is discussed in \cite{WNGuerreroBerrondo}.

The family of evolution operators  verifies the groupoid property:
\begin{eqnarray}
 U(t_2,t_1) U(t_1,t_0) &=& U(t_2,t_0)  \nonumber \\
 U(t_1,t_0)&=&U(t_0,t_1)^{-1} \label{groupoid}
\end{eqnarray}
$\forall t_0,t_1,t_2\in {\cal I}$. \modif{ See \cite{WNGuerreroBerrondo} for a discussion of the groupoid property in the context of the Wei-Norman factorization.}

Only in the time-independent case $A(t)=A\,,\forall t\in {\cal I}$ we find that $U(t,t_0)=U(t-t_0)$ and $U(t)$ defines
a 1-parameter group:
\begin{eqnarray}
 U(t_1+t_2)=U(t_1)U(t_2) \label{group}
\end{eqnarray}
with $A$ as infinitesimal generator.

Eq. \eqref{PVIop} is more general than Eq.  \eqref{PVIvec} and has the advantage that the initial condition, ${U}(t_0,t_0)={I}_V$, is fixed.
For simplicity of notation we shall simply write $U(t)\equiv U(t,t_0)$ and take $t_0=0$ most of the cases, but we should keep in mind the dependence of $U(t)$ on the initial time $t_0$.
 
In what follows we shall assume that $A(t)$ can be written as a finite sum:
\begin{equation}
 A(t)=\sum_{i=1}^n \alpha_i(t)X_i\,,\qquad \forall t\in{\cal I}
\end{equation}
where ${\cal B}=\{X_i,\,i=1,\ldots,n=\textrm{dim}{\cal G}\}$ forms a basis of a Lie algebra ${\cal G}$, realized as endomorphisms of $V$, with commutation relations given by:
\begin{equation}
 [X_i,X_j]=\sum_{k=1}^n C_{ijk} X_k \label{structureconstants}
\end{equation}
and  $C_{ijk}$ are known as the structure constants of the Lie algebra ${\cal G}$ and the basis ${\cal B}$.

Then the WN theorem \cite{WeiNorman1,WeiNorman2} states that the evolution operator $U(t)$ can be factorized as a product of exponentials and that this construction is valid in a neiboughood of $t_0$:


\begin{thm} \label{WNTheorem} \textit{(Local WN Theorem)}
 There is a neiboughood ${\cal J}\subset {\cal I}$ of $t_0$, where we can write $U(t)$ as:
 \begin{equation}
  U(t)= e^{g_1(t)X_1}e^{g_2(t)X_2}\cdots e^{g_n(t)X_n}\,,\qquad g_i(t_0)=0\,,\,i=1,\ldots,n \label{WNfactor}
 \end{equation}
with the functions $g_i(t)\,,i=1,\ldots,n$ satisfing a sistem of non-linear first-order ordinary diffrential equations, depending solely on the functions $\alpha_i(t)\,,i=1,\ldots,n$ and on the structure constants  of the Lie algebra ${\cal G}$.

\end{thm}

There are some important cases where the factorization is valid for all $t\in {\cal I}$ (global theorem), namely both for solvable algebras and  the case of $2\times 2$ real matrices \cite{WeiNorman2}.

Since we shall need it in subsequent sections, let us review the local theorem and skectch its proof (for the finite-dimensional case).

\noindent\textit{Proof:}

Solving for $A(t)$ in Eq. (\ref{PVIop}) and substituting $U(t)$ according to  Eq. (\ref{WNfactor}), we arrive to the following equation:

\begin{eqnarray}
\sum_{i=1}^n \alpha_i(t)X_i &=& \sum_{j=1}^n g_j'(t)\left(\prod_{k=1}^{j-1} e^{g_k(t)X_k}\right)X_j \left(\prod_{k=j-1}^1 e^{-g_k(t)X_k}\right) = \sum_{j=1}^n g_j'(t)
 \prod_{k=1}^{j-1} e^{g_k(t)ad_{X_k}}X_j   \nn \\
 &=& \sum_{j,k}^n g_j'(t) \xi_{jk}(g_1(t),g_2(t),\ldots,g_{j-1}(t))X_k \,\label{alphagp}
 \end{eqnarray}
 \noindent where $ad_X(Y)=[X,Y]$ is the \textit{adjoint} endomorphism of the Lie algebra
${\cal G}$.  The functions $\xi_{jl}$ are analytic functions of the $g_i$ depending solely on the structure constants $C_{ijk}$.

Since the $X_i$ constitute a basis of the Lie algebra, equating coefficients we arrive to the matricial equation $\alpha(t)=\Xi(g) g'(t)$:

\begin{equation}
 \left(
 \begin{array}{c}
        \alpha_1(t)\\
        \vdots \\
        \alpha_n(t)
\end{array}\right)
=
\left(
\begin{array}{ccc}
        \xi_{11}& \cdots & \xi_{n1}\\
        \vdots & \vdots &\vdots \\
        \xi_{1n}& \cdots & \xi_{nn}
\end{array}\right)
\left(
\begin{array}{c}
        g_1'(t)\\
        \vdots \\
        g_n'(t)
\end{array}\right) \,,
\left(
\begin{array}{c}
        g_1(t_0)\\
        \vdots \\
        g_n(t_0)
\end{array}\right) =
\left(
\begin{array}{c}
        0\\
        \vdots \\
        0
\end{array}\right) \label{alphagpmatrix}
\end{equation}

The matrix $\Xi(g)$ verifies $\Xi(g(t_0))=I_V$. By the analyticity of the functions $\xi_{jk}$, there is a neiboughood ${\cal K}$ of $t_0$ such that det$(\Xi(g))\neq 0$.

This allows to solve for $g'(t)$ in ${\cal K}$:
\begin{equation}
 g'(t) = \Xi(g){\,}^{-1} \alpha(t)\,,\qquad t\in {\cal K}
\end{equation}

Then we arrive to the IVP in ${\cal K}$:
\begin{equation}
 g'(t) = \Xi(g)^{-1} \alpha(t)\,,\qquad g(t_0)=0 \label{WeiNormanPVI}
\end{equation}

 Picard theorem tells us that there exists a neighborhood ${\cal J}\subset{\cal I}\cap{\cal K}$ where a unique solution of the IVP exists.$\blacksquare$

 Thus, assuming that the original IVP in Eq. (\ref{PVIvec}) has a unique solution on ${\cal I}=\mathbb{R}$, we only need to find the solutions of the equation det$(\Xi(g(t))= 0$ to obtain the intervals where the factorization (\ref{WNfactor}) exists.
 
 When det$(\Xi(g(t^*))= 0$, for some $t^*$, applying Rouche-Frobenius theorem, the linear system (\ref{alphagpmatrix}) can be undetermined or inconsistent. In the first case the $\alpha_i(t^*)$ do not determine completely the $g_i'(t^*)$, resulting in the loss of one or more degrees of freedom. In the second case, there are no values of $g_i'(t^*)$ consistent with the given $\alpha_i(t^*)$, resulting in the impossibility of describing the propagator $U(t^*)$ with the chosen factorization. When the group is $SO(3)$ parameterized by Euler angles, this phenomenon is usually called \textit{gimbal lock} \cite{Shuster}. For the case of the group $SU(1,1)$, in the context of the time dependent (and damped) harmonic oscillator, the values $t^*$ are associated with caustics, and are related to the Maslov correction and Morse index (see \cite{WNGuerreroBerrondo} for an explanation and references).

 The reason why for some values $t=t^*$ the propagator $U(t^*)$ can not be written as the factorization given in eq. (\ref{WNfactor}) should be understood in terms of local charts on the Lie group $G$.
 Eq. (\ref{WNfactor}) is nothing other than a local coordinate chart of the Lie group around the identity, known as canonical coordinates of the second kind on the group $G$ \cite{Varadarajan}.
 
 For a given choice of Lie algebra basis and an order of the basis elements, we have a different local coordinate chart. $|Det(\Xi(g))|$ can be interpreted as the Haar measure on $G$ in this local chart, but also as the Jacobian determinant of the transformation $O\subset \mathbb{R}^n \rightarrow G$ associated with eq. (\ref{WNfactor}). When $Det(\Xi(g))=0$ the Jacobian is not invertible and the transformation fails to be a homeomorphfism.
 
In this paper we shall apply the Wei-Norman factorization method to the Frenet-Serret equations \cite{Frenet,Serret} as an IVP on the Lie group $SO(3)$ of rotations.
The Frenet-Serret equations have many applications not only in differential geometry (curves on $\mathbb{R}^2$, $\mathbb{R}^3$ and other curved spaces, as well as their generalization  to curves in arbitrary dimensions), but also in other areas like  control of vehicles (bicycles, cars, airplanes) \cite{SeligRational}  or articulated robots \cite{Seligbook}, where
the group $SO(3)$ is enlarged to the Lie group $SE(3)$ of Euclidean transformations on $\mathbb{R}^3$.

 An important concept in control theory of vehicles or robots is that of spatial forces (forces and torques), also denoted wrenches \cite{Seligbook}. They are elements of the dual of the Lie algebra  $so(3)$ (or $se(3)$), usually called co-algebra and denoted by $so(3)^*$ (or $se(3)^*$). They are functionals on the Lie algebra, characterized by the evaluation map $W(\omega)=\langle W,\omega\rangle = Tr(W^\dag \omega),\, \omega\in so(3),\, W\in so(3)^*$ (where the Hilbert-Schmidt inner product is used).
 
 We can translate this concept to a general Lie group. Given the basis ${\cal B}$, defining the dual basis ${\cal B}^*=\{W_i,i=1,\ldots,n\}$, with $\langle W_i,X_j\rangle =\delta_{ij}$, the wrenches associated with the transformation $U(t)$ are given by:
 \begin{equation}
  \alpha_i(t)=\langle W_i,A(t) \rangle=\langle W_i,\frac{d U(t)}{dt} U(t)^{-1} \rangle = \sum_{j=1}^n g_j'(t)\xi_{ji}(g_1(t),\ldots,g_{j-1}(t))
  \label{wrenches}
 \end{equation}
 which is obtained by appling the duals $W_i$ to Eq. (\ref{alphagp}) or (\ref{alphagpmatrix}). Note that the coefficients $\alpha_i(t)$ play the role of forces, and they determine the derivatives $g_j'(t)$.
 
 Suppose that for $i_0\in I_0\subset \{1,2,\ldots,n\}$ we have $\alpha_{i_0}(t)\equiv 0$ for all $t$. Associated with these null wrenches we have the constraints:
 \begin{equation}
  \sum_{j=1}^n g_j'(t)\xi_{ji_0}(g_1(t),\ldots,g_{j-1}(t))=0\,,\qquad i_0\in I_0 
  \label{nullwrenches}
 \end{equation}
 
 These equations are automatically satisfied for solutions of (\ref{WeiNormanPVI}), but have the particularity that they do not depend on any wrenches $\alpha_i(t)$, i.e. they are \textit{structural equations}.
 They can be useful in designing trajectories and motion planning, using an \textit{anzat} for the functions $g_i(t)$ depending on parameters which are fixed imposing the initial and final conditions $y(t_0)$ and $y(t_f)$  and imposing the null wrenches (\ref{nullwrenches}). The resulting functions $g_i(t)$ are inserted in Eq. (\ref{alphagpmatrix}) and the remaining forces $\alpha_i(t), \, i\not\in I_0$ are derived to produce the desired motion \cite{SeligRational}.

 Particularizing eq. (\ref{WeiNormanPVI}) to the initial time $t_0$, we have the following results:
 
  \begin{cor}
   The initial velocities $g'(t_0)$ are given by:
   \begin{equation}
    g'_i(t_0) = \alpha_i(t_0) \,,\qquad i=1,2,\ldots, n
   \end{equation}
In particular, for the null wrenches:
   \begin{equation}
    g'_{i_0}(t_0) = 0 \,,\qquad i_0 \in I_0
   \end{equation}
   \label{CorollaryInitial}
  \end{cor}

 According to this result, at first-order in time, the evolution of the function $g_i(t)$ is governed by the corresponding wrench $\alpha_i(t)$. If the Lie algebra ${\cal G}$ were Abelian, this would be true for all orders. The non-commutativity of the Lie algebra causes the dependence on $g_i(t)$ on other wrenches $\alpha_j(t)$ for higher-orders in time.
 
 Another consequence of Corollary \ref{CorollaryInitial} is that, if any of the wrenches verifies $\alpha_i(t_0)\neq 0$, then 
 $g_i'(t_0)\neq 0$  and therefore there exist a neighbourhood ${\cal I}_0$ of $t_0$ where $g_i(t)$ is invertible, for all indices $i$ satisfying $g_i'(t_0)\neq 0$. Using this in eq. (\ref{nullwrenches}), we obtain:
 
 \begin{cor}
 For a fixed $j_0\not\in I_0$, we have:
   \begin{equation}
    \sum_{j=1,\,j\neq j_0}^n \frac{dg_j}{dg_{j_0}} \xi_{ji_0}(g_1,\ldots,g_{j-1})  = -\xi_{j_0i_0}(g_1,\ldots,g_{j_0-1})
   \end{equation}
   \label{CorollaryFunctionalRelation}
  \end{cor}

This Corollary provides a functional relationship between the variables $g_j$ and the variables $j_0\not\in I_0$, valid in a neighbourhood ${\cal I}_0$ of $t_0$. This relationship is independent of the wrenches. We will discuss an application of this corollary in the context of the Wei-Norman method for $SO(3)$ applied to solve the Frenet-Serret equations.

\section{Frenet-Serret equations }
\label{FSE}

The Frenet-Serret equations \cite{Frenet,Serret} describe the evolution of a particle moving (with non-zero velocity) in a differentiable trajectory in $\mathbb{R}^3$. From the geometrical point of view, the Frenet-Serret equations describe the properties of differentiable curves (with non-zero tangent vector), without relation to  any motion.

Consider a differentiable curve $\vec{r}(t),\,t\in\mathbb{R}$ in $\mathbb{R}^3$, with $\vec{r}{\,}'(t)\neq 0,\,\forall t\in\mathbb{R}$. Let us introduce a reparametrization of the curve in terms of the arclength $s$, defined as:
\begin{equation}
 s(t)=\int_{t_0}^t \|\vec{r}{\,}'(t')\|dt'
\end{equation}
where $t_0$ is the initial time, such that $s(t_0)=0$.

Since $\vec{r}{\,}'(t)\neq 0$, $s(t)$ is strictly increasing with $t$, and therefore the relation can be inverted and write $t(s)$, allowing to reparametrize the curve as $\vec{r}(s)=\vec{r}(t(s))$. For arclength parametrized curves the velocity is always one, $\|\frac{d \vec{r}(s)}{ds}\| =1$.

Therefore we can define the \textit{unitary tangent} vector
\begin{equation}
 \vec{T}(s)= \frac{d \vec{r}(s)}{ds}
\end{equation}

Since $\vec{T}(s)$ is unitary, its derivative with respect to $s$, if non-zero,  is orthogonal to itself. Suppose now that $\frac{d \vec{T}(s)}{ds}\neq 0$, meaning that the curve is nondegenate, i.e. the curvature is never zero. Then we can define the \textit{unitary normal} vector as:
\begin{equation}
 \vec{N}(s)= \frac{\frac{d \vec{T}(s)}{ds}}{\|\frac{d \vec{T}(s)}{ds}\|}
\end{equation}

The \textit{unitary binormal} vector is simply defined as $\vec{B}(s)=\vec{T}\times \vec{N}$.
By definition, the three vector $(\vec{T}(s),\vec{N}(s),\vec{B}(s))$ are unitary and orthogonal, defining, therefore, and orthonormal frame in $\mathbb{R}^3$ for every $s\in\mathbb{R}$. Since this frame is attached to the point $\vec{r}(s)$ of a curve, for $s\in\mathbb{R}$, it is usally called \textit{moving frame}.

The  Ferret-Serret equations describe how this moving frame changes along the curve:
\begin{eqnarray}
 \frac{d \vec{T}(s)}{ds}&=& \kappa(s)\vec{N}(s)  \nn\\
 \frac{d \vec{N}(s)}{ds}&=&  -\kappa(s) \vec{T}(s) + \tau(s) \vec{B}(s) \label{Eq:FSE} \\
  \frac{d \vec{B}(s)}{ds}&=&  -\tau(s)\vec{N}(s)  \nn
\end{eqnarray}

The first equation is simply the definition of $\vec{N}$, with $\kappa(s)=\|\frac{d \vec{T}(s)}{ds}\|$. The last equation is a consequence of the first equation and the orthonormality of the frame, with $\tau(s)=-\|\frac{d \vec{B}(s)}{ds}\|$. The second equation is a consequence of the other two equations. Note that if $\tau(s)\neq 0$ the curve is  non-flat, otherwise we can restric to a plane and ignore the third vector $\vec{B}(s)$.

If we consider the IVP formed by Eq. (\ref{Eq:FSE}) and the initial condition  $\vec{T}(s_0)=\vec{T}_0,\vec{N}(s_0)=\vec{N}_0,\vec{B}(s_0)=\vec{B}_0$, we see that these equations belong to the general form of Eq. (\ref{PVIvec}), formulated in the vector space $V=\mathbb{R}^9$.

Although the Frenet-Serret equations do not imply that $(\vec{T}(s),\vec{N}(s),\vec{B}(s))$ form an orthonormal frame in $\mathbb{R}^3$, as the following derived equations show (we have omitted the dependence in $s$):
\begin{eqnarray}
 \frac{d \vec{T}^2}{ds} &=& 2 \kappa  \vec{T}\cdot \vec{N}        \nn \\
 \frac{d \vec{N}^2}{ds} &=& -2  \kappa  \vec{T}\cdot \vec{N} + 2 \tau \vec{N}\cdot \vec{B}    \nn \\
 \frac{d \vec{B}^2}{ds} &=&  -2 \tau   \vec{N}\cdot \vec{B}     \\
 \frac{d (\vec{T}\cdot \vec{N})}{ds} &=&  \kappa(\vec{N}^2-\vec{T}^2) +\tau \vec{T}\cdot \vec{B} \nn \\
 \frac{d (\vec{T}\cdot \vec{B})}{ds}  &=& \kappa \vec{N}\cdot\vec{B} -\tau \vec{T}\cdot\vec{N}      \nn \\
 \frac{d (\vec{N}\cdot \vec{B})}{ds} &=& \tau(\vec{B}^2-\vec{N}^2) -\kappa\vec{T}\cdot\vec{B}     \nn
\end{eqnarray}
it is clear that if the initial conditions $(\vec{T}_0,\vec{N}_0,\vec{B}_0)$ form an orthonormal frame, then  $(\vec{T}(s),\vec{N}(s),\vec{B}(s))$ constitute an orthonormal frame $\forall s\in \mathbb{R}$. In fact it is enough that $(\vec{T}_0,\vec{N}_0,\vec{B}_0)$ are equal-norm orthogonal vectors, then $(\vec{T}(s),\vec{N}(s),\vec{B}(s))$ will be orthogonal and conserve the initial norm). The preserved properties $\vec{N}^2-\vec{T}^2=0$, $\vec{B}^2-\vec{N}^2=0$, $\vec{T}\cdot \vec{N}=0$, $\vec{N}\cdot \vec{B}=0$ and $\vec{T}\cdot \vec{B}=0$ are called invariant relations \cite{LeviCivita}, and restrict in a consistent way the dynamics to a submanifold of equal-norm orthogonal vectors.

Restricting the Frenet-Serret equations to the (nonlinear) manifold of orthonormal frames, and defining the orthogonal matrices 
\begin{equation}
 Q(s)=\left(
\begin{array}{c}
\vec{T}(s)^t \\
\vec{N}(s)^t \\
\vec{B}(s)^t
\end{array}
\right)
\end{equation}
i.e. with the vectors   $\vec{T}(s),\vec{N}(s),\vec{B}(s)$ written as rows, the Frenet-Serret equations can be rewritten as:
\begin{equation}
\frac{d Q(s)}{ds} = A(s) Q(s)\,,\qquad Q(s_0)=Q_0\in SO(3) \label{FSE-SO3}
\end{equation}
with
\begin{equation}
 A(s)=\left(
\begin{array}{ccc}
 0 & \kappa (s) & 0 \\
 -\kappa (s) & 0 & \tau (s) \\
 0 & -\tau (s) & 0 \\
\end{array}
\right)
\end{equation}

$A(s)$ is an antisymmetric matrix (a consequence of the orthonormality of the frame), and therefore is an element of
the Lie algebra of $SO(3)$, which will be denoted $so(3)$. An impotant feature of the Frenet-Serret frame, as opposed to other orthonormal frames, is that $A(s)$ depends on only two functions, $\kappa(s)$ and $\tau(s)$, instead of the three free parameters of a general antysymmetric matrix (see below for other frames sharing this feature). This means that the Frenet-Serret equations on $SO(3)$ has one null wrench, which can be used for motion planning.

Clearly, equations (\ref{FSE-SO3}) are written in the general form of Eq. (\ref{PVIvec}), but now formulated in the (nonlinear) manifold $SO(3)$. We restric to $SO(3)$ instead of considering the case $O(3)$ since the non-zero velocity condition implies that the orientation of the frame is preserved. Again, if we do not impose the initial condition $Q_0\in SO(3)$ the equations would be formulated in the vector space of $3\times 3$ matrices (isomorphic to $\mathbb{R}^9$).

It is worth commenting that the Frenet-Serret frame is not the only possibility to build an orthonormal moving frame, and other moving frames exist in the literature \cite{Bishop}  (with their corresponding equations) which possesses even better properties  in the sense that there is no need for the regularity condition \cite{Bishop} and the resulting frame is more well-behaved \cite{Mansfield}, preserving the feature of only two functions in the antysimmetric matrix $A(s)$.

Also, the Frenet-Serret equations can be easily generalized to higher-dimensions or to other manifolds like spheres \cite{Struik} or space forms \cite{FERRANDEZ20061666}.

\section{Wei-Norman method applied to Ferret-Serret equations on $SO(3)$}

According to Eqn. (\ref{FSE-SO3}), the Frenet-Serret equations can be formulated as a linear first-order PVI on $SO(3)$.  Let us write $Q(s)=U(s,s_0)Q(s_0)$, with $U(s,s_0)\in SO(3)$ the propagator\footnote{Note that since the original PVI in Eq. (\ref{FSE-SO3}) is already in matrix form,  we could refrain from introducing the propagator. We will do so in order to stay within the Wei-Norman framework.},
 which rotates the frame $Q(s)$ along the curve, where $Q(s)$ is interpreted as a frame (which is itself a rotation matrix). The equation for the propagator $U(s,s_0)$ is simply:
\begin{equation}
 \dparcial{U(s,s_0)}{s}={A}(s){U}(s,s_0)\,,\qquad {U}(s_0,s_0)={I}_3 \label{PVIopSO(3)}
 \end{equation}
 
 Once we have written our original problem into this form we shall apply the Wei-Norman method to reformulate the matrix differential equation (\ref{PVIopSO(3)}) (involving 9 unknowns) into a simpler set of coupled linear differential equations involving 3 unknowns.
 
 For that purpose we have to introduce an ordered basis in the Lie algebra $so(3)$, as indicated in Sec. \ref{WN}. We shall work with two different bases, the standard basis $\{L_x,L_y,L_z\}$ where $L_i$ generates rotations around the $i$th axis:
 \begin{equation}
  L_x=\left(\begin{array}{ccc}
       0&0&0\\
       0&0&-1 \\
       0&1&0
      \end{array}\right), \quad
L_y=\left(\begin{array}{ccc}
       0&0&1\\
       0&0&0 \\
       -1&0&0
      \end{array}\right), \quad  
L_z=  \left(\begin{array}{ccc}
       0&-1&0\\
       1&0&0 \\
       0&0&0
      \end{array}\right), \quad 
 \end{equation}
 and the complex basis $\{L_+,L_-,L_x\}$, with $L_\pm = L_z \pm i L_y$.
 
 Although in the Wei-Norman factorization method all elements in the Lie algebra basis are included, cf. Eq.  (\ref{WNfactor}), for semisimple Lie algebras there can be repeated generators (omitting others), and in fact this is common in the case of $so(3)$ where Euler angles are used \cite{Shuster}. We shall consider here both type of factorizations, Euler rotations (of types z-x-z, x-y-x, y-z-y, z-y-z, x-z-x, y-x-y)  and Cardan-Tait-Bryan rotations (of types x-y-z, y-z-x, z-x-y, x-z-y, z-y-x, y-x-z) \cite{GeneralizedEuler}, this last one corresponding to the usual Wei-Norman factorization of $SO(3)$.
 
 As in Sec. \ref{WN}, to simplify the notation we shall write $U(s)$ instead of $U(s,s_0)$, but we should keep in mind the dependence of the propagator on the inicial contition $s_0$ (which in most of the cases it will be taken as $s_0= 0$). 
 
 \subsection{Basis $\{L_x,L_y,L_z\}$ }
 
 For the Cardan-Tait-Bryan angles  we shall consider  the order z-x-y, since it produces the simplest equations, together with the order z-y-x, where the propagator is written as:
 \begin{equation}
  U(s)=e^{\psi(s)L_z}e^{\theta(s)L_x}e^{\phi(s)L_y}\,,\label{WNfactorSO3-real}
 \end{equation}
with $\psi,\theta,\phi\,: {\cal J}\rightarrow [0,2\pi)$. The Wei-Norman equations (see \cite{WeiNorman1,WeiNorman2} or Theorem \ref{WNTheorem}) associated with this factorization are:
 \begin{eqnarray}
 \psi '(s)&=&-\kappa (s)-\tau (s) \tan \theta (s) \sin \psi (s) \,,\qquad \psi(s_0)=0\nn \\
 \theta '(s)&=& -\tau (s)\cos\psi (s) \,,\qquad\qquad\qquad\qquad\,\,\; \theta(s_0)=0\label{WNeqSO3-real}\\
 \phi '(s)&=&\tau (s) \sec \theta (s) \sin \psi (s) \,,\qquad\qquad\qquad \phi(s_0)=0 \nn
  \end{eqnarray}
  
 Note that with this factorization the equations are particularly simple when $\tau(s)=0$ (planar curves), reducing to $\psi '(s)=-\kappa (s)$ which is inmediatly integrable by a quadrature. For other cases the equations are highly non-linear and coupled and therefore difficult or impossible to solve.
 
 The Wei-Norman equations for the order z-y-x are obtained from Eq. (\ref{WNeqSO3-real}) through the change $\psi(s)\rightarrow \psi(s)+\frac{\pi}{2},\,\theta(s)\rightarrow \theta(s) +\pi$.  The Wei-Norman equation for the orders x-z-y and x-y-z and  are obtained from  z-x-y  and z-y-x by the changes $\kappa(s)\leftrightarrow \tau(s),\,\theta(s)\rightarrow \theta(s)+\pi$ and $\kappa(s)\leftrightarrow \tau(s),\,\theta(s)\rightarrow -\theta(s)$, respectively.
 
 For the  Euler angles we shall consider  the order z-y-z (again it produces the simplest equations, together with the order z-x-z), where the propagator is written as:
 \begin{equation}
  U(s)=e^{\psi(s)L_z}e^{\theta(s)L_y}e^{\phi(s)L_z}\,.\label{WNfactorSO3-real-Euler}
 \end{equation}
with $\psi,\theta,\phi\,: {\cal J}\rightarrow [0,2\pi)$.  The Wei-Norman equations associated with this factorization are:
 \begin{eqnarray}
 \psi '(s)&=&-\kappa (s)+\tau (s) \cot \theta (s) \cos \psi (s) \nn \,,\qquad\quad \psi(s_0)=0\\
 \theta '(s)&=& \tau (s)\sin\psi (s)\,,\qquad\qquad\qquad\qquad\qquad\,\,\; \theta(s_0)=0 \label{WNeqSO3-real-Euler}\\
 \phi '(s)&=&-\tau (s) \csc \theta (s) \cos \psi (s)\,,\qquad\qquad\qquad\,\, \phi(s_0)=0 \nn
  \end{eqnarray}
  
 These equations are equivalent to the previous case with the change $\psi(s)\rightarrow \psi(s)+\frac{\pi}{2},\theta(s)\rightarrow\theta(s)+\frac{\pi}{2}$.
 
 It is well-known that for semisimple Lie groups like {$SO(3)$ the validity of these differential equations is local \cite{WeiNorman1,WeiNorman2}. In fact, we can state the following Corollary, whose proof is immediate particularizing Theorem \ref{WNTheorem} to this case (see also \cite{WeiNorman1,WeiNorman2}).
 
 \begin{cor}
 The solution of the first-order PVI (\ref{WNeqSO3-real} )(resp. (\ref{WNeqSO3-real-Euler})) is only valid in intervals not containing the values of $s$ such that $\cos\theta(s)= 0$ (resp. $\sin\theta(s)= 0$).
 \end{cor}

 Although the Eqs. (\ref{WNeqSO3-real}) (resp. (\ref{WNeqSO3-real-Euler})) are ill-defined for $\cos\theta(s)= 0$ (resp. $\sin\theta(s)= 0$), which are known as ''caustic points`` in the context of the group $SU(1,1)$ \cite{WNGuerreroBerrondo}, the resulting propagator is well-defined. This is a general issue of the Wei-Norman method, i.e. even if the factorization is not valid, the solution of the initial linear first-order PVI (the Frenet-Serret equations in this case) exist for all values of $t$ if the curvature and torsion are continuous functions. In general we can compute it performing carefully a limit \cite{WNGuerreroBerrondo}. For instance, for the propagator (\ref{WNfactorSO3-real}) and $\theta(s^*)=(2k+1)\frac{\pi}{2}\,,\,k=0,1$, we have:
 \begin{equation}
 \lim_{s\rightarrow s^*} U(s)=\left(
\begin{array}{ccc}
 \cos (\psi (s^*)+(-1)^k\phi (s^*)) & 0 & \sin (\psi (s^*)+(-1)^k\phi (s^*)) \\
 (-1)^k\sin (\psi (s^*)+(-1)^k\phi (s^*)) & 0 & -(-1)^k\cos (\psi (s^*)+(-1)^k\phi (s^*)) \\
 0 & (-1)^k & 0 \\
\end{array} \label{Us-GimbalLock}
\right)
\end{equation}

Even though  $\psi(s^*)$ and $\phi(s^*)$ are angles representing rotation wit respect to independent axes, $U(s^*)$ in Eq. (\ref{Us-GimbalLock}) represents a rotation with respect to a single axis by the sum or difference of  $\psi(s^*)$ and $\phi(s^*)$, thus a degree of freedom is lost. In addition, under the propagator (\ref{Us-GimbalLock}) the vector $\vec{N}$ transforms to $(-1)^k \vec{N}$. This lost of one degree of freedom can produce, in certain circustances, that $\ldots$

This is a drawback of this kind of factorizations, associated with singularities in the chosen coordinates, usually known in the literature as \textit{Gimbal lock} \cite{Shuster}, and it always affects  the central rotation ($\theta$ in both cases, although the rotation axes are different).

From the topological point of view, the Gimbal lock is a consequence of the fact that there is no local homeomorphism from $\mathbb{T}^3$ (parametrized by  the 3 angles $(\psi,\theta,\phi)$) into $SO(3)$, which coincides with $\mathbb{R}P^3$. The WN factorization in terms of Euler or Cardan angles can be seen a local chart in the group SO(3), more precisely a \textit{local chart of the second kind}.

Further insight can be gained in the gimbal lock if we apply the Rouché-Frobenious theorem to the linear system in eq. (\ref{alphagpmatrix}). For the Cardan angles parametrization with the order z-x-y we have:

 \begin{thm}(Rouch\'e-Frobenious applied to Wei-Norman on $SO(3)$)\\
 The Wei-Norman equations are:
 \begin{itemize}
  \item[(i)] Compatible (rank$\Xi(s)=3$) if $\cos\theta(s)\neq 0$.
  \item[(ii)] Underdetermined ($\cos\theta(s^*)= 0$, rank$\Xi(s^*)=2$ and rank of enlarged matrix $= 2$ ) if $\sin(\psi(s^*))\tau(s^*) = 0$.
  \item[(iii)] Incompatible ($\cos\theta(s^*)= 0$, rank$\Xi(s^*)=2$ and rank of enlarged matrix $= 3$ ) if $\sin(\psi(s^*))\tau(s^*) \neq 0$.
 \end{itemize}
 \end{thm}

Thus, in the undetermined case, the gimbal lock is not dangerous, it is only a coordinate singularity that can be avoided performing carefully the limit as in eq. (\ref{Us-GimbalLock}) (and fixing a particular solution of the undetermined system). However, in the incompatible case, there is a real lock since the applied wrenches (the torsion in the considered coordinartes) has an orthogonal component to the distribution generated by  $\Xi(s^*)$. This is related to the lock experienceb by a three-axes gimbal in a certain configuration of the rings \cite{Shuster}.

In Engineering, Gimbal lock appears when using gimbals to monitor and control the orientation of aircrafts (or spaceships) or in controling robotic arms, where it is known as \textit{wrist flip}. A way to avoid the Gimbal Lock is to use gimbals with four axes (instead of three) and controlling manually the fourth axis to keep the others orthogonal to avoid their lock \cite{Shuster}.

In computer science (particularly in computer graphics), when rotations are computed, it is more advangeous to use quaternions, which are defined on $\mathbb{S}^3$ (a double covering of $SO(3)$). They define a map from $\mathbb{S}^3$ to  $SO(3)$ which are true local homeomorphism and therefore they do not suffer from Gimbal lock. However, care must be taken of the fact that a rotation can be associated with two distinct quaternions, and the extra overhead of complex arithmetics. In most cases, Tait-Bryan or Euler rotations are used and when a Gimbal lock is approaching a switch to quaternions is made and then back to Tait-Bryan or Euler rotations once the Gimbal lock is avoided.

In our case, a way to avoid the Gimbal lock is to change to a different local chart (for instance changing the factorization from  (\ref{WNfactorSO3-real}) to (\ref{WNfactorSO3-real-Euler}) and vice versa) whenever a Gimbal lock is approaching.

\subsection{Basis $\{L_+,L_-,L_x\}$ }

The reason to use a complex basis is to \textit{uncouple} the highly coupled set of differential equations obtained with the real basis $\{L_x,L_y,L_z\}$ in the previous subsection. This is essentially equivalent to the complex coordinate change introduced by Lie and Darboux (see pag. 36 in \cite{Struik}) which uncouples the equations into a Riccati equation and two quadratures (once the Riccati equation has been solved). We shall easily recover  these equations using the complex basis $\{L_+,L_-,L_x\}$ in a suitable order. The factorization with the order $(L_+,L_x,L_-)$:
 \begin{equation}
  U(s)=e^{g_+(s)L_+}e^{g_x L_x}e^{g_-L_-}\,,\label{WNfactorSO3-complex}
 \end{equation}
leads to the equations:
\begin{eqnarray}
g_+{}'(s)&=&-\frac{1}{2} \left(1+g_+(s)^2\right) \kappa(s)-i g_+(s)\tau (s) \,,\qquad g_+(s_0)=0\nn\\
g_x{}'(s)&=&-\tau(s)+i g_+(s) \kappa(s)\,,\qquad\qquad\qquad\qquad\;\; g_x(s_0)=0\label{WNeqSO3-complex}\\
g_-{}'(s)&=&-\frac{1}{2}\kappa (s)e^{i g_x(s)}\,,\qquad\qquad\qquad\qquad\qquad \quad g_-(s_0)=0\nn
\end{eqnarray}
which is an uncoupled Riccati equation and the other two equations can be solved by quadratures.
The order $(L_-,L_x,L_+)$ leads to the same equatons interchanging $g_+ \leftrightarrow g_-$. The orders $(L_\pm,L_\mp,L_x)$ lead to an uncoupled  Riccati equation and the other two equations are slightly more complex although solvable by quadratures. The orders $(L_x,L_\pm,L_\mp)$ lead to coupled equations.

The price we have to pay for using a complex basis is that the unknowns $(g_\pm(x),g_x(s))$ are complex, therefore there are six (real) degrees of freedom. Even more,  $U(s)$ is not, explictly,  a real orthogonal matrix. Clearly, there must be some contraints on the functions  $(g_\pm(x),g_x(s))$ for  $U(s)$ to be a real matrix. Imposing Im$(U(s))=0$, we obtain the following constraints:
\begin{eqnarray}
 |g_+(s)|&=&|g_-(s)|\equiv \rho(s) \nn \\
 \textrm{Im}(g_x(s))&=&-\log \frac{|g_+(s)|(1+|g_-(s)|^2)}{|g_-(s)|}=-\log(1+\rho(s)^2) \label{RealityConstraint} \\
 \textrm{Re}(g_x(s))&=& 2\pi k+\phi_+(s)+\phi_-(s)\nn
\end{eqnarray}
where $g_\pm(s)=|g_\pm(s)|e^{i\phi_\pm(s)}$, $k\in \mathbb{Z}$, $\rho(s)\geq 0$, $\phi_\pm(s)\in[0,2\pi)$.

The system of complex equations (\ref{WNeqSO3-complex}) reduce to the real system:
\begin{eqnarray}
 \rho '(s)&=& -\frac{1}{2}  \cos \phi_+(s)\,(1+\rho(s)^2) \kappa (s)\,,\qquad\qquad\quad \rho(s_0)=0\nn\\
 \phi_+{}'(s) &=& -\frac{ \left(\rho (s)^2-1\right) \sin\phi_+(s)}{2 \rho (s)}\kappa (s)-\tau (s) \,,\qquad \phi_+(s_0)=0 \\
 \phi_-{}'(s)&=& -\frac{ \left(\rho (s)^2+1\right) \sin \phi_+(s)}{2 \rho (s)}\kappa (s)\,,\qquad\qquad\quad \phi_-(s_0)=0 \nn
\end{eqnarray}

The initial condition for the angles $\phi_+$ and $\phi_-$ are arbitrary since they are not defined for $\rho=0$. These equations are again highly non-linear and coupled, therefore we shall not work with them.

We shall work with the equations (\ref{WNeqSO3-complex}), keeping in mind the constraints (\ref{RealityConstraint}), which can be written in the form of invariant relations \cite{LeviCivita}.

The Riccati equation, being difficult or impossible to solve, it is one of the best known equations and there are many ways of obtaining solutions (see for instance \cite{Struik}).
It is interesting its relation with the second order linear equations, established through the change:
\begin{equation}
 g_+(s) = \frac{2 u'(s)}{\kappa (s) u(s)}
\end{equation}
with $u(s)$ satisfying the complex second-order PVI:
\begin{equation}
 u''(s)-\left(\frac{\kappa '(s)}{\kappa (s)}-i \tau (s)\right) u'(s)+\frac{1}{4} \kappa (s)^2 u(s)\,,\qquad u(s_0)=1\,,\qquad u'(s_0)=0 \label{2orderEq}
\end{equation}

This equation can be interpreted as a complex  harmonic oscillator of $s$-dependent frequency $\omega(s)=\frac{1}{2} \kappa(s)$ and complex  $s$-dependent damping term.
The damping term can be eliminated if we perform the new change:
\begin{equation}
 u(s)=\sqrt{\frac{\kappa (s)}{\kappa(s_0)}}  e^{-\frac{1}{2} i \int_{s_0}^s \tau (t)\, dt}  v(s)\,,
\end{equation}
with $v(s)$ satisfying:
\begin{eqnarray}
 v''(s) &+& \left(\frac{\kappa (s)^2}{4} +\frac{\tau (s)^2}{4} -\frac{3 \kappa '(s)^2}{4 \kappa (s)^2}+\frac{\kappa ''(s)}{2 \kappa (s)}+ \frac{i}{2 \kappa (s)}W(\tau,\kappa)(s) \right)v(s) =0 \nn \\
 v(s_0)&=&1 \label{2orderEq-real} \\
v'(s_0)&=& -\frac{\kappa'(s_0)}{2 \kappa (s_0)} +i \frac{\tau (s_0)}{2} \nn
\end{eqnarray}
where $W(f,g)(s)=f(s)g'(s)-f'(s)g(s)$ the Wronskian of $f$ and $g$. If the Wronskian of $\kappa$ and $\tau$ is not zero, $v(s)$ satisfies the equation of a harmonic oscillator with a complex frequency. Otherwise the frequency is real. Note that if the Wronskian is zero $\tau(s)/\kappa(s)$ is constant, and therefore the curve is of constant slope or Lancret helix \cite{libroAngel}.


The solutions of the Eq. (\ref{WNeqSO3-complex}) can be expressed in terms of the solutions of the Eq. (\ref{2orderEq}) and (\ref{2orderEq-real}):
\begin{eqnarray}
 g_+(s)&=& \frac{2 u'(s)}{\kappa (s) u(s)} = \frac{2 v'(s)}{\kappa (s) v(s)}+\frac{\kappa '(s)}{\kappa (s)^2} -\frac{i \tau (s)}{\kappa (s)} \nn\\
 g_x(s) &=& -\int_{s_0}^s \tau (t)dt+2 i \log u(s) = 2 i \log \left( v(s)\sqrt{\frac{\kappa(s)}{\kappa(s_0)}}\right) \\
 g_-(s) &=& -\frac{1}{2}\int_{s_0}^s \frac{\kappa (t) e^{-i \int_{s_0}^t \tau (w) \, dw}}{u(t)^2}  dt = -\frac{\kappa(s_0)}{2}\int_{s_0}^s \frac{1}{v(t)^2}  dt\nn
\end{eqnarray}


\subsubsection{Real parametrization: Stereographic projection coordinates}

For the complex factorization (\ref{WNfactorSO3-complex}), let us introduce real coordinates that implement the constraints (\ref{RealityConstraint}):
\begin{eqnarray}
g_+(s) &=& \tan \frac{\theta (s)}{2} e^{i \phi_+(s)} \nn\\
g_x(s) &=& \phi_+(s)+\phi_-(s)+2 i \log \left|\cos\frac{\theta (s)}{2}\right|
\label{SO3-stereographic}\\
g_-(s)&=& \tan \frac{\theta (s)}{2} e^{i \phi_-(s)} \nn
\end{eqnarray}
where $\theta(s)\in[0,\pi)$.
With this coordinates, $g_+(s)$ can be interpreted as the stereographic projection coordinates (from the North pole) of a point moving on the Riemann sphere with lattitude angle $\theta(s)-\frac{\pi}{2}$ and azimutal angle $\phi_+(s)$. The interpretation of $g_-(s)$ is the same, but changing only the azimutal angle (i.e. the two points  are always at the same lattitude).

In terms of the real coordinates defined by Eq. (\ref{SO3-stereographic}), the Wei-Norman equations are:
\begin{eqnarray}
 \theta '(s) &=&-\kappa (s)\cos \phi_+(s) \nn  \\
 \phi_+'(s) &=& \kappa (s) \cot \theta (s) \sin \phi_+(s)-\tau (s) \label{WNeqSO3-stereographic} \\
 \phi_-'(s) &=& -\kappa (s) \csc \theta (s) \sin \phi_+(s) \nn
\end{eqnarray}

It coincides with the Euler factorization with the order x-z-x with the identification 
$\theta(s)\rightarrow \theta(s),\phi_+(s)\rightarrow \psi(s), \phi_-(s)\rightarrow \phi(s)$.

For the complex basis, the Wei-Norman Local Theorem states:

\begin{cor}
 The solution of the system of first-order PVI (\ref{WNeqSO3-complex}) is only valid in intervals not containing the values of $s$ such that the inverse image of $g_+(s)$ (and $g_-(s)$) under the stereographic projection does not pass through the north pole of the Riemann sphere.
\end{cor}
\textit{Proof.} Using Theorem \ref{WNTheorem} we compute $\Xi(g)$, resulting in det$(\Xi(g))=-2 i e^{-i g_x(s)}$. This determinant vanishes when $| g_+(s)| = | g_-(s)| =\rho(s)\rightarrow \infty$, which corresponds to $\theta(s)=\pi$, i.e., to the north pole of the Riemann sphere.$\blacksquare$

Since the initial values $g_+(s_0)=g_-(s_0)=0$ correspond to the south pole of the Riemann sphere, the factorization (\ref{WNfactorSO3-complex}) fails for the values $s^*$ when the unit tangent verifies $\vec{T}(s^*)=-\vec{T}(s_0)$. In fact, we have for the propagator:%
\begin{equation}
\lim_{s\rightarrow s^*} U(s)=\left(
\begin{array}{ccc}
 -1 & 0 & 0 \\
 0 & -\cos (\phi_-(s^*)-\phi_+(s^*)) & \sin (\phi_-(s^*)-\phi_+(s^*)) \\
 0 & \sin (\phi_-(s^*)-\phi_+(s^*)) & \cos (\phi_-(s^*)-\phi_+(s^*)) \\
\end{array}
\right)
\end{equation}

Note that a similar phenomenom to Gimbal Lock appears here, since $U(s^*)$ depends on the difference $\phi_-(s^*)-\phi_+(s^*)$, and therefore one degree of freedom is lost.

\section{Conclusions}
\label{sec:conclusions}

We have applied the Wei--Norman factorization method to the Frenet--Serret
equations, regarded as a linear initial value problem on the rotation group
$SO(3)$. The method replaces the matrix equation for the propagator, with
nine unknowns, by a system of three nonlinear first-order equations for the
coordinates of the propagator in a chart of the second kind on the rotation group. The form of this
system depends strongly on the choice and ordering of the Lie algebra basis.

With the real basis $\{L_x,L_y,L_z\}$, the Cardan and
Euler  factorizations immediately integrable for plane curves, but are highly coupled otherwise.
Their singular points, $\cos\theta(s)=0$ and $\sin\theta(s)=0$ respectively,
are the gimbal lock of the corresponding angle parameterization. We have
shown that at these points the Wei--Norman equations break down but the
propagator has a well-defined limit, which depends only on the sum or the
difference of the outer angles, so that one degree of freedom is lost. Since these singularities belong to the chart
rather than to the solution, they can be avoided by switching between
factorizations as a singular point is approached. 

With the complex basis $\{L_+,L_-,L_x\}$ and the order $(L_+,L_x,L_-)$, the
Wei--Norman equations uncouple into a Riccati equation for $g_+(s)$ followed
by two quadratures. This recovers, within a Lie-algebraic framework, the
classical reduction of Lie and Darboux \cite{Struik}. The price is that
the exponents are complex, so we have derived the constraints that ensure
the propagator is a real orthogonal matrix. The Riccati equation can be
linearized into a complex harmonic oscillator with frequency $\kappa(s)/2$,
whose damping term can be removed. The resulting effective frequency is
real exactly when the Wronskian $W(\tau,\kappa)$ vanishes, that is, for
Lancret helices \cite{libroAngel}. A real parameterization of the
complex factorization identifies $g_\pm(s)$ with stereographic coordinates
on the Riemann sphere and turns out to coincide with one of the Euler 
factorizations. In this geometric picture the factorization fails exactly
when the unit tangent is reversed, $\vec T(s)=-\vec T(s_0)$, and the limiting
propagator again shows the loss of one degree of freedom.


Several questions remain open.  The structural
constraints associated with the null wrench of the Frenet--Serret frame
could be exploited for motion planning of vehicles and robots, in the spirit
of \cite{SeligNeedle,SeligRational}, and the whole analysis could be extended to
the Euclidean group $SE(3)$. Other moving frames with only two independent
coefficients, such as the rotation minimizing (Bishop) frame
\cite{Bishop,Mansfield}, can be treated with the same method.
Finally, the generalization to higher dimensions and to curves in space
forms \cite{FERRANDEZ20061666} would require the Wei--Norman method for $SO(n)$ and
for the pseudo-orthogonal groups, where the interplay between the choice of
charts and the singularities of the factorization becomes richer.

\section*{Acknowledgements}

We thank the support of the Spanish MICINN through the project PID2022-138144NB-I00.

\bibliographystyle{unsrt}
\bibliography{References}

\end{document}